\documentstyle[12pt]{article}

\def\C{C\!\!\!\!I}
\def\dna{\downarrow}

\def\O{{\cal O}} 
\def\PP{{\bf P}} 
\def\T{{\cal T}}
\def\Z{Z\!\!\!Z}

\def\hra{\hookrightarrow}
\def\id{\mathop{\rm id}\nolimits} 
\def\im{\mathop{\rm im}\nolimits} 
\def\ker{\mathop{\rm ker}\nolimits}
\let\lra\longrightarrow
\let\pa\partial
\let\ov\overline

\newtheorem{theorem}{Theorem}[subsection]
\newtheorem{proposition}[theorem]{Proposition}
\newtheorem{lemma}[theorem]{Lemma}

\def\rem{\refstepcounter{theorem}\paragraph{Remark \thetheorem}}

\def\proof{\paragraph{Proof}}
\def\example{\refstepcounter{theorem}\paragraph{Example \thetheorem}}

\makeatletter \def\l@section{\@dottedtocline{1}{0em}{1.2em}} \makeatother

\begin{document}

\title{Characteristic Classes for $GO(2n,\C)$}
\author{Yogish I. Holla and Nitin Nitsure}

\date{}

\maketitle

School of Mathematics, Tata Institute of Fundamental Research,
Homi Bhabha Road, Mumbai 400 005, India.

\tableofcontents

\section{Introduction}

The complex Lie group $GO(n)$ is by definition the 
closed subgroup of $GL(n)$ consisting of all
matrices $g$ such that ${^t}gg$ is a scalar matrix
$\lambda I$ for some $\lambda \in \C^*$.
(We write simply $GL(n)$, $SO(n)$, $O(n)$, etc. for
the complex Lie groups $GL(n,\C)$, $SO(n,\C)$, $O(n,\C)$, etc.) 

The group $GO(1)$ is just $\C^*$, with 
classifying space $B\C^*=\PP^{\infty}_{\C}$, whose
cohomology ring $H^*(B\C^*;\Z/(2))$ is the
polynomial ring $\Z/(2)[\lambda]$
where $\lambda\in H^2(B\C^*;\Z/(2))$ 
is the Euler class.
For an odd number $2n+1\ge 3$, 
the group $GO(2n+1)$ is isomorphic to
the direct product $\C^*\times SO(2n+1)$. 
Hence $BGO(2n+1)$ is homotopic to the direct product 
$B\C^*\times BSO(2n+1)$, with cohomology ring 
the polynomial ring
$\Z/(2)[\lambda,w_2,w_3,\ldots,w_{2n+1}]$ where 
for $2\le i\le 2n+1$, the elements $w_i \in H^i(BSO(2n+1);\Z/(2))$
are the Stiefel-Whitney classes.

In this paper, we consider the even case $GO(2n)$.
The main result is an 
explicit determination, in terms of generators and
relations, of the singular cohomology ring $H^*(BGO(2n);\Z/(2))$.
This is the Theorem \ref{cohomology} below.

An outline of the argument is as follows. 

To each action of the group $\C^*$ on any space $X$, we
functorially associate a certain derivation 
$s : H^*(X)\to H^*(X)$ on
the cohomology ring of $X$, which is graded 
of degree $-1$, with square zero 
(see Lemma \ref{derivation via action} below). 
In terms of the action
$\mu : \C^* \times X\to X$ and the projection 
$p : \C^* \times X\to X$, it is given by the formula
$$\mu^* - p^* = \eta\otimes s$$ 
where $\eta$ is the positive generator of $H^1(\C^*)$. 

When $X$, together with its given $\C^*$-action, 
is the total space of a principal
$\C^*$-bundle $\pi : X\to Y$ over some base $Y$, recall that we have a
long exact Gysin sequence 
$$\cdots \stackrel{\lambda}{\to}  
H^i(Y) \stackrel{\pi^*}{\to} H^i(X) 
\stackrel{d}{\to} H^{i-1}(Y)
\stackrel{\lambda}{\to} H^{i+1}(Y) 
\stackrel{\pi^*}{\to}  \cdots$$
We show (Lemma \ref{derivation via Gysin}) that in this case, 
the derivation $s$ on $H^*(X)$ equals the composite
$$s = \pi^*\circ d$$

As $O(2n)\subset GO(2n)$ is normal with quotient $\C^*$,
$BO(2n)$ can be regarded as the total space of
a principal $\C^*$-bundle $BO(2n)\to BGO(2n)$.
The resulting action of $\C^*$ on $BO(2n)$ gives 
a derivation $s$ on the ring 
$H^*(BO(2n)) = \Z/(2)[w_1,\ldots,w_{2n}]$
with $\mu^* - p^* = \eta\otimes s$.
The last equality enables us to write the following
expression for $s$ (see Lemma \ref{s is s'} below)
$$s = \sum_{i=1}^n w_{2i-1} {\pa \over \pa w_{2i}}$$
Determining the invariant subring $\ker(s)$ of $H^*(BO(2n))$
thus becomes a purely commutative algebraic problem, which
we solve by using the technique of regular sequences and Koszul complex,
and express the invariant subring $B=\ker(s)$ in terms
of generators and relations (see Theorem \ref{y2k} below).

As explained above, in terms of the Gysin sequence
{\small
$$\cdots \stackrel{\lambda}{\to} H^i(BGO(2n)) 
\stackrel{\pi^*}{\to} H^i(BO(2n)) 
\stackrel{d}{\to} H^{i-1}(BGO(2n)) 
\stackrel{\lambda}{\to} H^{i+1}(BGO(2n)) 
\stackrel{\pi^*}{\to} \cdots$$
}%end of small
the derivation $s$ equals the composite $\pi^*\circ d$.  
This fact, together with the explicit knowledge
of the kernel and image of $s$ (in terms of the ring $B$), 
allows us to `solve' the above long exact Gysin sequence,
to obtain the main Theorem \ref{cohomology}, which gives
generators and relations for the ring $H^*(BGO(2n))$.

Note that $H^*(BGL(2n))$ is the polynomial ring
$\Z/(2)[\ov{c}_1,\ldots,\ov{c}_{2n}]$ in the variables
$\ov{c}_i$ which are the images (mod $2$) of the Chern classes
$c_i \in H^{2i}(BGL(2n);\Z)$. The natural
inclusion $GO(2n)\hra GL(2n)$ induces a ring homomorphism 
$H^*(BGL(2n)) \to H^*(BGO(2n))$, and we determine the images of the 
$\ov{c}_i$ in $H^*(BGO(2n))$ in terms of our generators of 
$H^*(BGO(2n))$ (see Proposition \ref{odd Chern classes} below). 

The paper is arranged as follows. The purely algebraic 
problem of determining the kernel ring of the derivation 
$s$ is solved in Section 2.
In Section 3, we describe the derivation on the cohomology of a space  
associated to a $\C^*$-action, and connect this with the 
Gysin sequence for principal $\C^*$-bundles.
The above material is applied to 
determine the cohomology ring of $BGO(2n)$ in Section 4. 
In Section 5, we determine the ring homomorphism 
$H^*(BGL(2n)) \to H^*(BGO(2n))$ which is induced by the natural
inclusion $GO(2n)\hra GL(2n)$. (This section is new in this version,
and does not appear in [H-N].)

\section{Derivation on a polynomial algebra}

Let $k$ denote the field $\Z/(2)$.
For any integer $n\ge 1$, consider the 
polynomial ring $C= k[w_1,\ldots,w_{2n}]$ where
$w_1,\ldots,w_{2n}$ are independent 
variables (the ring $C$ will stand for $H^*(BO(2n))$,
with $w_i$ the Stiefel-Whitney classes, when 
we apply this to topology).

In what follows, we will treat the odd variables 
$w_{2i-1}$ differently from the even variables $w_{2i}$, so 
for the sake of clarity we introduce the notation
$$u_i = w_{2i-1}~~\mbox{and } v_i = w_{2i}~~\mbox{for }
1\le i \le n.$$
On the ring $C=k[u_1,\ldots,u_n, v_1,\ldots,v_n]$, 
we introduce the derivation
$$s = \sum_{i=1}^n u_i {\partial \over \partial v_i} : C\to C$$
Note that $s$ is determined by the conditions
$s(u_i) = 0$ and $s(v_i) = u_i$ for 
$1\le i \le n$. As $2=0$ in $k$, it follows from the above that
$$s\circ s = 0$$
We are interested in the subring $B = \ker(s)\subset C$.
Let $A\subset B \subset C$ be the subring
$A = k[u_1,\ldots,u_n, v_1^2,\ldots,v_n^2]$.
For any subset $T=\{ p_1, \ldots, p_r\} \subset \{ 1,\ldots,n\}$
of cardinality $r\ge 1$, consider the monomial
$v_T=v_{p_1}\cdots v_{p_r}$. We put $v_T=1$ when $T$ is empty.
It is clear that $C$ is a free $A$-module of
rank $2^n$ over $A$, with $A$-module basis formed by the
$v_T$. Hence as the ring $A$ is noetherian,
the submodule $B\subset C$ is also finite over $A$.
In particular, this proves that $B$ is a finite type $k$-algebra.

We are now ready to state the main result of this
section, which is a description of the ring $B$ by
finitely many generators and relations.

\begin{theorem}\label{y2k} 
{\bf (I) (Generators) } Let the derivation
$s : C \to C$ be defined by 
$s= \sum_{i=1}^n u_i {\partial \over \partial v_i}$.
Consider the subring $B=\ker(s) \subset C$.  
We have $s^2 =0$, and $B = \ker(s) = \im(s) +  A$.
Consequently, 
$B$ is generated as an $A$-module by $1$ together with the 
$2^n-n-1$ other elements $s(v_T)$, where $T$ ranges over
all subsets $T \subset  \{ 1, \ldots, n\}$
of cardinality $\ge 2$, and $v_T = \prod_{i\in T}v_i$.

{\bf (II) (Relations) } Let $A[c_T]$ be the polynomial ring over 
$A$ in the $2^n-n-1$ 
algebraically independent variables $c_T$, indexed by 
all subsets $T \subset \{ 1,\ldots,n\}$ of cardinality $\ge 2$.
Let $N$ be the ideal in $A[c_T]$ generated by 

(1) the $2^n -n(n-1)/2 -n-1$ elements of the form
$$\sum_{i\in T} u_i\,c_{T-\{ i\}}$$
where $T \subset \{ 1,\ldots,n\}$ is a subset of cardinality $\ge 3$, and

(2) the $n(n-1)/2$ elements of the form
$$(c_{\{ i,j \} })^2 + u_i^2v_j^2 + u_j^2v_i^2$$
where $\{ i,j\}\subset \{ 1,\ldots,n\}$ is a subset of cardinality $2$, and

(3) the $(2^n -n -1)^2 - n(n-1)/2$ elements 
(not necessarily distinct) of the form
$$c_T\,c_U + \sum_{p\in T} \prod_{q \in T \cap U -\{ p\} }  
\, u_p \, v_q^2 \, c_{(T - \{ p\}) \Delta U}$$
where $T\ne U$ in case both $T$ and $U$ have cardinality $2$,
and $\Delta$ denotes the symmetric difference of sets 
$X \Delta Y = (X-Y)\cup (Y-X)$.

Then we have an isomorphism of $A$-algebras
$A[c_T] /N \to B$, mapping $c_T \mapsto s(v_T)$.
\end{theorem}

{\bf Proof } of (I) (Generators) : 
As already seen, $s\circ s=0$, hence 
we have $\im(s) \subset \ker(s) = B$. 
As $A\subset B$, we have the inclusion 
$A + \im(s) \subset B$. We next prove that in fact 
$B = A + \im(s)$. The problem translates into the exactness of a 
certain Koszul resolution, as follows. An exposition of the
elementary commutative algebra used below can be found, 
for example, in \S 16 of Matsumura [M].

Let $M$ be a free $A$-module of rank $n$, with 
basis $e_1,\ldots,e_n$. Consider the $A$-linear map 
$u: M \to A : e_i \mapsto u_i$. 
For any subset $T=\{ p_1, \ldots, p_r\} \subset \{ 1,\ldots,n\}$
of cardinality $r$ where $0\le r\le n$, let
$e_T = e_{p_1} \wedge \ldots \wedge e_{p_r} \in \bigwedge^rM$.
These $e_T$ form a free $A$-module basis of $\bigwedge^r M$.
Recall that the Koszul complex for $u$ is the complex
$$0 \lra \bigwedge^n M \stackrel{d_n}{\lra} \bigwedge^{n-1} M 
\stackrel{d_{n-1}}{\lra} \cdots \stackrel{d_3}{\lra} 
\bigwedge^2 M  \stackrel{d_2}{\lra}  M  \stackrel{d_1}{\lra} A \lra 0$$
where for $1\le r \le n$, the differential 
$d_r : \bigwedge^r M \to \bigwedge^{r-1} M$ is defined by putting
\begin{eqnarray}d_r(e_T) = \sum_{j\in T} u_j e_{T - \{ j\}}
\end{eqnarray} 
(note that $T$ is non-empty as $r\ge 1$). 
As $u_1, \ldots, u_n$ is a regular sequence 
in $A$, the Koszul complex is exact except in degree $0$
(in fact the Koszul complex gives rise to a projective resolution
of the $A$-module $A/\im(d_1) = A/(u_1,\ldots,u_n)$,  
but we do not need this).

Let $\bigwedge(M) = \oplus_{0\le r \le n}\bigwedge^rM$, which we regard as
a free $A$-module of rank $2^n$ with basis $e_T$.
(The algebra structure of $\bigwedge(M)$ is not relevant to us.)  
The graded module $\bigwedge(M)$ comes with an $A$-linear endomorphism
$d$ of degree $-1$, where by definition $d$ is $d_r$ on
the graded piece $\bigwedge^rM$. Hence we have the equality
$\ker(d) = \im(d) + A$ as submodules of $\bigwedge(M)$.

Now consider the $A$-linear isomorphism $\psi: \bigwedge(M) \to C$
of free $A$-modules under which $e_T \mapsto v_T$.
From the expression $s = \sum_i u_i(\pa /\pa v_i)$, it follows that
\begin{eqnarray}s(v_T)  = \sum_{j\in T} u_j v_{T - \{ j\}}
\end{eqnarray}
when $T$ is non-empty.  
By comparing (1) and (2) we see that 
$\psi: \bigwedge(M) \to C$ takes the Koszul differential $d$
to the differential $s$, hence $\ker(s) = \im(s) + A$. 

As $\im(s)$ is generated by all the elements $s(v_T)$ where
$T$ is non-empty, the $A$-module $B$ is generated by
$1$ together with the $2^n-1$ other elements $s(v_T)$, where $T$ ranges over
all non-empty subsets $T \subset  \{ 1, \ldots, n\}$.

This completes the proof of part (I) (Generators) of the Theorem \ref{y2k}.

{\bf Proof } of Theorem \ref{y2k}.(II) (Relations) : 
As $s\circ s =0$, the equation (2) 
gives $\sum_{i\in T} u_i s(v_{T - \{ i\} }) =0$.
This proves that
the relation II.(1) is satisfied by the assignment $c_T \mapsto s(v_T)$.

In particular, when $T = \{ i,j\}$ has cardinality $2$, we get
$s(v_{\{ i,j \} })  =  u_iv_j + u_jv_i$, and hence
$s(v_{\{ i,j \} })^2  =  u_i^2v_j^2 + u_j^2v_i^2$. This proves that
the relation II.(2) is satisfied by the assignment $c_T \mapsto s(v_T)$.

Note that if $X$ and $Y$ are subsets of $\{ 1,\ldots, n\}$,
then we have
\begin{eqnarray}v_X v_Y = \prod_{i\in X\cap Y} v_i^2 \,v_{X\Delta Y}
\end{eqnarray}
where $X\Delta Y$ denotes the symmetric difference
$(X-Y)\cup (Y-X)$.

From the equations (2) and (3) it follows by a straight-forward
calculation that 
$$s(v_T)\,s(v_U)  =  
\sum_{p\in T} \prod_{q \in T \cap U -\{ p\} }  
\, u_p \, v_q^2 \, s(v_{(T - \{ p\}) \Delta U})$$
where $T\ne U$ in case both $T$ and $U$ have cardinality $2$.
This proves that the relation II.(3) is also satisfied by the 
assignment $c_T \mapsto s(v_T)$.

Hence the $A$-algebra homomorphism $A[c_T] \to B$ defined by sending
$c_T \mapsto s(v_T)$ kills the ideal $N$, hence defines
an $A$-algebra homomorphism $\eta : A[c_T]/N \to B$. 
By part (I) of the theorem, which we have already proved, 
the $A$-algebra $B$ is generated by the $s(v_T)$, 
hence $\eta$ is surjective.  

To prove that $\eta$ is injective, we use the following lemma.

\begin{lemma} 
Let $F \subset C$ be the $A$-submodule generated
by all $v_T$ where $T\subset \{ 1,\ldots,n\}$ has cardinality
$\ge 1$. 
Let $E \subset B$ be the $A$-submodule generated
by all $s(v_T)$ where $T\subset \{ 1,\ldots,n\}$ has cardinality
$\ge 2$. Then we have $E\subset F$, and $F\cap A =0$, and 
$B = A \oplus E$, as $A$-modules.
\end{lemma}

{\bf Proof } Note that for each $i$, 
the $v_i$-degree of each non-zero monomial in 
any element of $A$ is even. On the other hand, if $|T| \ge 1$
(where $|T|$ denotes the cardinality of $T$), 
then the monomial $v_T$ has $v_i$-degree $1$ for
for all $i\in T$. Hence in any $A$-linear combination 
$w=\sum_{|T|\ge 1} a_T\,v_T$, each monomial
has odd $v_i$-degree for some $i$. It follows that
$A\cap F =0$. 

From the equality 
$s(v_T) =  \sum_{i\in T} u_i v_{T - \{ i\} }$, 
it follows that if $|T| \ge 2$ then
$s(v_T) \in F$, so $E \subset F$. It follows that
$A\cap E =0$. By part (I) of Theorem \ref{y2k}, we 
already know that $A + E =B$,
so the lemma follows. \hfill$\Box$

\medskip

We now prove that $\eta$ is injective. For this, we define a homomorphism
of $A$-modules $\theta \,: \,C \to A[c_T]/N$ by
$\theta(1) = 0$, $\theta(v_i) = u_i$, and 
$\theta(v_T) = c_T$ for $|T| \ge 2$.
With this definition of $\theta$, note that the composite
$C \stackrel{\theta}{\to} A[c_T]/N
\stackrel{\eta}{\to} B \hra C$
equals $s : C \to C$. Moreover, note that 
$\theta (E) = 0$, 
where $E \subset C$ is the $A$-submodule in the previous lemma.
To see the above, we must show that $\theta(s(v_T)) = 0$ for all $T$.
For $|T|\ge 3$, this follows from the relations (1) in part (II) of
Theorem \ref{y2k}, while for $|T|\le 2$, it follows from 
$\theta(v_i)=u_i$ and $\theta(1) =0$. 

Now suppose $x\in A[c_T]/N$ with $\eta(x) = 0$. 
From the generators (3) of $N$, it follows that any element of
$x\in A[c_T]/N$ can be written in the form
$$x = a_0 + \sum_{|T|\ge 2} a_T \,c_T~~~
\mbox{where }a_0,\, a_T \in A$$
Hence by definition of $\eta$, we have
$\eta(x) = a_0 +  \sum_{|T| \ge 2} a_T\,s(v_T)$.
As $\eta(x)=0$, and as $A\cap E =0$ by the above lemma, 
we have 
$$a_0 = \sum_{|T| \ge 2} a_T\,s(v_T) = 0$$
Hence we get
$$x = \sum_{|T|\ge 2} a_T \,c_T = \theta(x')~~\mbox{where }
x' = \sum_{|T|\ge 2} a_T \,v_T$$ 
As $s(x') = \sum_{|T| \ge 2} a_T\,s(v_T) = 0$, 
we have $x' \in B$. 
In the notation of the above lemma, $x'$ is in $F$, 
and as $B\cap F =E$, we see that $x' \in E$.
As $\theta(E)=0$, it
finally follows that $x = \theta(x') = 0$.

This proves $\eta$ is injective, and completes the proof
of the Theorem \ref{y2k}. \hfill$\Box$

\section{$\C^*$-actions, derivations, and Gysin boundary}

\subsection{The derivation associated to a $\C^*$-action}

Given an action $\mu : \C^*\times X \to X$ of the group $\C^*$ on a 
space $X$, consider the pullback
$\mu^* : H^*(X) \to H^*(\C^*\times X)$, where the cohomology 
is with arbitrary coefficients.
By K\"unneth formula, for each $n$ we get a homomorphism
$\mu^* : H^n(X) \to H^0(\C^*) \otimes H^n(X) 
\oplus H^1(\C^*)\otimes H^{n-1}(Y)$. 
As we have a section $X \to \C^* \times X : x\mapsto (1,x)$,
it follows that for any $\alpha$ in $H^n(X)$, we have
$\mu^*(\alpha) = 1\otimes  \alpha+ \eta\otimes \alpha'$
where $\eta$ is the positive generator of $H^1(\C^*)$, and
the element $\alpha' \in H^{n-1}$ is uniquely determined by $\alpha$.
On the other hand, under the projection $p: \C^*\times X \to X$ 
we have $p^*(\alpha) = 1\otimes  \alpha$.

We now define a linear operator $s : H^n(X) \to H^{n-1}(X)$ 
on $H^*(X)$ by the equality
$\mu^*(\alpha) = 1\otimes\alpha + \eta\otimes s(\alpha)$.

\example Let $X=\C^*$, and consider the
action $m : \C^*\times \C^* \to \C^* : (x,y)\mapsto xy$. 
Then we have
the basic equality
$$m^*(\eta) = 1\otimes \eta + \eta\otimes 1 $$ 
It follows that the corresponding operator $s$ on $H^*(\C^*) =
k[\eta]/(\eta^2)$ is given by $s(1) = 0$ and $s(\eta) =1$.

\begin{lemma}\label{derivation via action}
Given an action $\mu : \C^*\times X \to X$ of the group
$\C^*$ on a space $X$, the
map $s : H^*(X)\to H^*(X)$ 
on the singular cohomology ring of $X$ with arbitrary coefficients, 
defined by the equality 
$$\mu^* -p^* = \eta\otimes s$$ 
where $p: \C^*\times X \to X$ is the projection 
and $\eta$ is the positive generator of $H^1(\C^*)$, 
is a graded anti-derivation of degree
$(-1)$ on the ring $H^*(X)$, that is,  
for all $\alpha \in  H^i(X)$ and $\beta \in H^j(X)$, we have
$$s(\alpha  \beta) = 
s(\alpha)  \beta + (-1)^i\,
\alpha \,  s(\beta) \in H^{i+j-1}(X)$$
Moreover, the composite $s\circ s =0$.
\end{lemma}

{\bf Proof } The verification of the derivation property is
straight from the definition. We now prove the property $s\circ s =0$.
For this consider the commutative diagram
$$\begin{array}{ccc}
\C^*\times\C^*\times X&\stackrel{(m,\id_X)} {\lra}
&\C^*\times X \\
\scriptstyle{(\id_{\C^*},\mu)}\dna~~~~ &
& ~~~~\dna\scriptstyle{\mu} \\
\C^*\times X&\stackrel{\mu}{\lra}
& X 
\end{array}$$
where $m : \C^*\times \C^* \to \C^* : (x,y)\mapsto xy$. 
As $m^*(\eta) = 1\otimes \eta + \eta\otimes 1$, it follows 
from the definition of $s$ that for any $\alpha \in H^i(X)$,
\begin{eqnarray*}
(m,\id_X)^*\circ \mu^* (\alpha) & = & 
1\otimes 1\otimes \alpha + 
1\otimes \eta \otimes s(\alpha) 
+ \eta\otimes 1 \otimes s(\alpha),~\mbox{ and } \\
(\id_{\C^*},\mu)^*\circ \mu^* (\alpha) & = & 
1\otimes 1\otimes \alpha + 
1\otimes \eta \otimes s(\alpha) + 
\eta\otimes 1 \otimes s(\alpha) + 
\eta\otimes \eta\otimes s^2(\alpha).
\end{eqnarray*}
Comparing, we get $\eta\otimes \eta\otimes s^2(\alpha)=0$, 
which means $s^2(\alpha)=0$.
\hfill$\Box$

\subsection{Derivation in the case of principal $\C^*$-bundles}

We next express the above derivation $s$ in the case when $X$
is the total space of a principal 
$\C^*$-bundle $\pi_X : X\to Y$ over some base $Y$, 
in terms of the Gysin boundary map $d_X : H^*(X)\to H^*(Y)$.

\begin{lemma}\label{derivation via Gysin}  
Let $\pi_X : X\to Y$ be a principal 
$\C^*$-bundle. Then the graded anti-derivation $s_X$ 
on the singular cohomology ring $H^*(X)$ with arbitrary
coefficients, which is associated to the given $\C^*$-action on $X$ 
by Lemma \ref{derivation via action}, equals the composite 
$$s_X = \pi_X^*\circ d_X$$ 
where $d_X : H^*(X) \to H^*(Y)$ is the Gysin boundary map
of degree $(-1)$ and $\pi_X : H^*(Y) \to H^*(X)$ is induced by
the projection $\pi_X : X\to Y$.
\end{lemma}

{\bf Proof } Let $\mu_X : \C^* \times X \to X$  
be the action, and $p_X : \C^* \times X \to X$ be 
the projection. If the bundle $X$ is trivial, then the equality 
$\mu_X^* - p_X^* = \eta\otimes s_X$ is obvious from the definitions. 
So it remains to prove this equality for a non-trivial $X$.

For any map $f : B \to Y$, let $\pi_M : M \to B$
denote the pullback bundle. 
The two projections $\pi_M :M\to B$ and
$r: M \to X$, and the Gysin boundary map $d_M$, satisfy
$$\pi_M^*\circ f^* = r^*\circ \pi_X^*~~\mbox{and }
d_M \circ r^* = f^*\circ d_X$$
Here the second equality follows from the fact that the Euler class of
$X$ pulls back to the Euler class of $M$.
Hence the following diagram commutes
$$\begin{array}{ccc}
H^n(X) & \stackrel{\pi_X^*d_X}{\lra} & H^{n-1}(X) \\
\scriptstyle{r^*}\dna~~ & & ~~\dna \scriptstyle{r^*}\\
H^n(M) & \stackrel{\pi_M^*d_M}{\lra} & H^{n-1}(M)
\end{array}$$

For the projection $(\id_{\C^*},r) : \C^*\times M\to \C^*\times X$,
we similarly have
$$p_M^* r^* = (\id_{\C^*},r)^*p_X^* ~~\mbox{and }
\mu_M^* r^* = (\id_{\C^*},r)^*\mu_X^*$$ 

Now take $f: B\to Y$ to be $\pi_X : X \to Y$.
Then $M$ is trivial and so the lemma holds for $s_M$. 
The projection $r:M\to X$
has a tautological section $\Delta : X\to M$ which is 
the diagonal of $M = X\times_Y X$. 
(Under the canonical isomorphism 
$(\mu_X,p_X): \C^* \times X \to X\times_Y X$, 
the section $\Delta : X\to M$ becomes the map $u \mapsto (1,u)$.)

Hence given
$\alpha\in H^i(X)$ we have
\begin{eqnarray*}
\eta\otimes \pi_X^*d_X(\alpha) 
& = & \eta\otimes \Delta^*r^* \pi_X^*d_X(\alpha)
      ~~\mbox{as }r\circ\Delta =\id_{X}.\\
& = & \eta\otimes \Delta^*\pi_M^*d_Mr^*(\alpha)
        ~~\mbox{as }r^*\pi_X^*d_X = \pi_M^*d_M r^*.\\
& = &  (\id_{\C^*}, \Delta)^* ( \eta\otimes \pi_M^*d_M r^*\alpha) \\
& = &  (\id_{\C^*}, \Delta)^* (\mu_M^* -p_M^*)(r^*\alpha) 
       ~~\mbox{as $M$ is trivial.} \\
& = & (\id_{\C^*}, \Delta)^* (\id_{\C^*}, r)^*(\mu_X^* -p_X^*)(\alpha)\\
& = & (\mu_X^* -p_X^*)(\alpha)~~\mbox{as }
       (\id_{\C^*}, r) \circ(\id_{\C^*}, \Delta) = \id_{\C^*\times X}.\\
& = & \eta\otimes s_X(\alpha)~~\mbox{by definition of }s_X.
\end{eqnarray*}
This proves the lemma in the general case.  \hfill$\Box$

\section{Cohomology of $BGO(2n)$}

{\bf Note } From now onwards, in the rest of this paper, 
singular cohomology will be with coefficients $k=\Z/(2)$,
unless otherwise indicated.

\subsection{Principle $GO(n)$-bundles and reductions to $O(n)$}

\centerline{\bf Principal $GO(n)$-bundles and triples $(E,L,b)$}

Recall that principal $O(n)$-bundles $Q$ on a space $X$ are equivalent
to pairs $(E,q)$ where $E$ is the rank $n$
complex bundle on $X$ associated to $Q$ via the defining representation of
$O(n)$ on $\C^n$, and $q:E\otimes E \to \O_X$ is the everywhere 
non-degenerate symmetric bilinear form on $E$ with values
in the trivial line bundle $\O_X$ on $X$, such that $q$ is
induced by the standard quadratic form $\sum x_i^2$ on $\C^n$. 
The converse direction of
the equivalence is obtained via a Gram-Schmidt process,
applied locally.

By a similar argument, principal $GO(n)$-bundles $P$ 
on $X$ are equivalent to triples $(E,L,b)$ where $E$ is 
the vector bundle associated to $P$ via 
the defining representation of
$GO(n)$ on $\C^n$, and $b : E\otimes E \to L$ 
is the everywhere nondegenerate
symmetric bilinear form on $E$ induced by the standard
form $\sum x_i^2$ on $\C^n$, which now takes values in the
line bundle $L$ on $X$, associated to $P$ by the homomorphism
$\sigma : GO(n)\to \C^*$ defined by the equality
${^t}gg = \sigma(g)I$.

\example\label{eachL} Let $L$ be any line bundle on $X$.
On the rank $2$ vector bundle $F = L\oplus \O_X$, 
we define a nondegenerate symmetric bilinear form $b$
with values in $L$ by putting
$b: (x_1,x_2)\otimes (y_1,y_2) \mapsto x_1\otimes y_2 + y_1\otimes x_2$.
Now for any $n\ge 1$, let the triple  
$(F^{\oplus n},L,b)$ be the orthogonal direct sum of $n$ copies of 
$(F,L,b)$. This shows that given any space $X$ 
and any even integer $2n\ge 2$,
there exists 
some nondegenerate symmetric bilinear triple $(E,L,b)$ on $X$ 
of rank $2n$, where $L$ is a given line bundle.

\medskip

\centerline{\bf Reductions to $O(n)$}

Given a principle $GO(n)$-bundle $P$ on $X$, let $(E,L,b)$ be the
corresponding triple, and let $L_o = L-X$ (complement of zero section).
As the sequence 
$$1\to O(n) \to GO(n) \stackrel{\sigma}{\to} \C^* \to 1$$
is exact, 
$L_o\to X$ is the associated $GO(n)/O(n)$-bundle to $P$, and
so reductions of structure group of $P$ from $GO(n)$ to $O(n)$ are
the same as global sections (trivializations) of $L$.
This can be described in purely linear terms, by saying that
given an isomorphism $v: L \stackrel{\sim}{\to} \O_X$, we get a
pair $(E, v\circ b : E\otimes E \to \O_X)$ from the triple
$(E,L,b)$.

Let $u:L \stackrel{\sim}{\to} \O_X$ and 
$v: L \stackrel{\sim}{\to} \O_X$ be two such
reductions. Then the two $O(n)$-bundles $(E, u\circ b)$ and $(E, v\circ b)$ 
are not necessarily isomorphic. In particular, the 
two sets of Stiefel-Whitney classes 
$w_i(E, u\circ b)$ and $w_i(E, v\circ b)$ 
need not coincide, but are related as follows. 

Given any rank $n$ vector bundle $E$ together
with an $\O_X$-valued nondegenerate bilinear form 
$q : E\otimes E \to \O_X$,
consider the new bilinear form $yq : E\otimes E \to \O_X$
where $y: X \to \C^*$ is a nowhere vanishing function. 
Let $(y) \in H^1(X,\Z/(2))$ be the pull-back of the 
generator $\eta \in H^1(\C^*,\Z/(2))$ (in other words, 
$(y)$ is the Kummer class of $y$). Note that $(y)^2=0$.

A simple calculation using the splitting principle shows the following.

\begin{lemma}\label{w of yq}
If $w_i(E, q)$ are
the Stiefel-Whitney classes of the $O(n)$-bundle $(E,q)$, then the 
Stiefel-Whitney classes of the $O(n)$-bundle $(E,yq)$ 
are given by the formula
$$w_i(E, yq) = w_i(E, q) + (n-i+1)(y)\cdot  w_{i-1}(E, q)$$
\end{lemma}

\rem It can be seen that the homomorphism $GO(2n)\to \{ \pm 1\} :
g \mapsto \sigma(g)^n/\det(g)$ is surjective, with connected 
kernel $GSO(2n)$. This implies that $\pi_0(GO(2n)) =\Z/(2)$,
and hence $\pi_1(BGO(2n))=\Z/(2)$.

\subsection{The natural derivation on the ring $H^*(BO(n))$}

Let $G$ be a Lie group and $H$ a closed subgroup. If $P \to BG$
is the universal bundle on the classifying space $BG$ of $G$, then
the quotient $P/H$ can be taken to be $BH$, so that
$BH \to BG$ is the bundle associated to $P$ by the action of
$G$ on $G/H$. In particular, taking $G=GO(n)$ and $H=O(n)$,
we see that $BO(n) \to BGO(n)$ is the fibration
$L_o \to BGO(n)$, where $L_o$ is the complement of the
zero section of the line bundle $L$ on $BGO(n)$ occuring in
the universal triple $(E,L,b)$ on $BGO(n)$.

We denote by $\pi : BO(n) \to BGO(n)$ the
projection and we denote by $\lambda\in H^2(BGO(n))$ 
the Euler class of $L$. This gives us the long exact Gysin sequence
{\small
$$\cdots \stackrel{\lambda}{\to}   
H^r(BGO(n)) \stackrel{\pi^*}{\to} 
H^r(BO(n)) \stackrel{d}{\to} 
H^{r-1}(BGO(n))\stackrel{\lambda}{\to} 
H^{r+1}(BGO(n))  \stackrel{\pi^*}{\to} \cdots $$
}%end of small

By Lemma \ref{derivation via Gysin}, the composite maps
$s = \pi^*\circ d : H^r(BO(n))\to H^{r-1}(BO(n))$
define a derivation $s$ on the graded ring $H^*(BO(n))$.
We now identify this derivation.

Recall that the singular cohomology ring $H^*(BO(n))$
with $k=\Z/(2)$ coefficients is the polynomial ring
$H^*(BO(n)) = k[w_1,\ldots,w_n]$ in 
the Stiefel-Whitney classes $w_i$.
We use the convention that $w_0=1$.

\begin{proposition}\label{s is s'}
Let $s : H^*(BO(n)) \to H^*(BO(n))$ be the composite
$s = \pi^*\circ d$ as above. Then $s$ is a derivation
of degree $(-1)$ on the graded ring $H^*(BO(n))$, with
$s\circ s =0$. 
In terms of the universal Stiefel-Whitney classes $w_i$,
we have $H^*(BO(n)) = k[w_1,\ldots,w_n]$, and the derivation
$s$ is given in terms of these generators by
$$s = \sum_{i=1}^n (n-i+1) w_{i-1}{\pa\over \pa w_i}
\, : \, w_i \mapsto (n-i+1)w_{i-1}$$
In particular, for $n=2m$ even, let $u_i = w_{2i-1}$ and $v_i =w_{2i}$,
where $1\le i\le m$ be the generators of the polynomial
ring $H^*(BO(2m))$. Then the derivation $s$ 
on the polynomial ring $k[u_1,\ldots u_m,v_1,\ldots,v_m]$
is given by $s= \sum u_i{\pa\over \pa v_i}$, 
with kernel ring explicitly 
given in terms of generators and relations by the Theorem \ref{y2k}. 
\end{proposition}

{\bf Proof } Let $(E,L,b)$ denote the universal triple on $BGO(n)$
and $\pi : L_o\to BGO(n)$ the projection, where
$L_o$ is the complement of the zero section of $L$.
The pullback $\pi^*(L)$ has a tautological trivialization 
$\tau : \pi^*(L) \stackrel{\sim}{\to} \O_{L_o}$,
which gives the universal pair $(\pi^*E, \tau\circ \pi^*q)$ 
on $L_o=BO(n)$. 
Let $f=\pi\circ p : \C^*\times L_o \to BGO(2n)$ 
be the composite of the projection $p : \C^*\times L_o \to L_o$ 
with $\pi : L_o\to BGO(n)$. 

Under the projection $p : \C^*\times L_o \to L_o$, the pair 
$(\pi^*E, \tau\circ \pi^*b)$ pulls back
to the pair $(f^*E, p^*(\tau\circ \pi^*q))$ on $\C^*\times L_o$.
Under the scalar multiplication 
$\mu : \C^*\times L_o \to L_o$, the pair 
$(\pi^*E, \tau\circ \pi^*q)$ pulls back
to the pair $(f^*E, yp^*(\tau\circ \pi^*b))$ on $\C^*\times L_o$,
where $y : \C^*\times L_o \to \C^*$ is the projection.
By Lemma \ref{w of yq}, we have 
$$w_i[\mu^*(\pi^*E, \tau\circ \pi^*b)] = 
w_i[p^*(\pi^*E, \tau\circ \pi^*b)] + 
(n-i+1) (y)\cdot  w_{i-1}[p^*(\pi^*E, \tau\circ \pi^*b)]$$
in the cohomology ring $H^*(\C^*\times L_o)$.

Note that the Stiefel-Whitney classes $w_i(\pi^*E, \tau\circ \pi^*b)$
are simply the universal classes $w_i$.
The class $(y)$ becomes $\eta\otimes 1$ under the K\"unneth
isomorphism. Hence the above formula reads 
$$\mu^*w_i = p^*w_i + (n-i+1)\eta\otimes w_{i-1}$$
Hence the proposition follows by lemmas \ref{derivation via action} and
\ref{derivation via Gysin}.
\hfill$\Box$

\subsection{Generators and relations for the ring $H^*(BGO(2n))$}

\centerline{\bf The elements $\lambda$, $a_{2i-1}$, 
$b_{4i}$ and $d_T$ of $H^*(BGO(2n))$}

Let $(E,L,b)$ be the universal triple on $BGO(2n)$.
Recall that we denote by $\lambda \in H^2(BGO(2n))$
the Euler class of $L$.
By Example \ref{eachL} applied to $X = \PP^{\infty}_{\C}$
with $L=\O_X(1)$ the universal line bundle on $X$,
together with the universal property of $BGO(2n)$, 
it follows that 
$\lambda^n\ne 0 ~~\mbox{for all } n\ge 1$.
For each $1\le j\le n$, we define elements 
$$a_{2j-1} = dw_{2j}\in H^{2j-1}(BGO(2n)).$$
Note that we therefore have 
$$\pi^*(a_{2j-1}) = s(w_{2j}) = w_{2j-1}\in H^{2j-1}(BO(2n)).$$
More generally, for any subset 
$T=\{ i_1,\ldots,i_r\} \subset \{ 1,\ldots,n\}$ 
of cardinality $r\ge 2$, let 
$v_T = w_{2i_1}\cdots w_{2i_r} \in H^{2\deg(T)}(BO(2n))$
where $\deg(T) = i_1 + \cdots + i_r$. We put
$$d_T = d(v_T) \in H^{2\deg(T) - 1}(BGO(2n)).$$
Next, let $E$ be the rank $2n$ complex vector bundle on 
$BGO(2n)$ which occurs in the universal triple $(E,L,b)$ on $BGO(2n)$.
For $1\le j\le n$, let $b_{4j} \in H^{4j}(BGO(2n))$ be 
the element 
$$b_{4j} = \ov{c}_{2j}(E)\in H^{4j}(BGO(2n))$$ 
which is the image of the $2j$ th Chern class 
$c_{2j}(E) \in H^{2j}(BGO(2n);\Z)$ of $E$ 
under the change of coefficients from $\Z$ to $\Z/(2)$. 
For any $m\ge 1$, consider the composite map 
$$H^*(BGL(m);\Z) \to H^*(BGL(m);\Z/(2)) \to H^*(BO(m);\Z/(2))$$ 
where the first map is induced by the change of coefficients 
$\Z \to \Z/(2)$ and the second map is induced by the natural inclusion
$O(m)\hra GL(m)$. The ring $H^*(BGL(m);\Z)$ is the polynomial ring
$\Z[c_1,\ldots,c_m]$ where $c_i\in H^{2i}(BGL(m);\Z)$ is the $i$ th
universal Chern class. The image of $c_i$ in 
$H^*(BO(m);\Z/(2))$ is known to be $w_i^2$, for all $1\le i\le m$. 
Hence the elements 
$b_{4j} \in H^{4j}(BGO(2n))$ satisfy the property that 
$$\pi^*(b_{4j}) = w_{2j}^2~~\mbox{for each } 1\le j\le n.$$  

In the ring $H^*(BGO(2n))$, we have 
$\lambda  a_{2i-1}=0$ for all $1\le i\le n$, and
$\lambda  d_T=0$ for every subset $T\subset \{ 1,\ldots,n\}$
of cardinality $|T| \ge 2$, as
follows from the fact that $\lambda \circ d =0$ 
in the Gysin sequence.

\medskip

\centerline{\bf The image of $\pi^*$}

We now come to a crucial lemma, one which allows us to write down the
image of the ring homomorphism $\pi^* : H^*(BGO(2n)) \to H^*(BO(2n))$.

\begin{lemma}\label{image of alpha}
We have the equality
$\ker(d) = \im(\pi^*) = \ker(s) \subset H^*(BO(2n))$, 
in other words, the sequence 
$H^i(BGO(2n)) \stackrel{\pi^*}{\to} 
H^i(BO(2n))\stackrel{s}{\to} H^{i-1}(BO(2n))$
is exact.
\end{lemma}

{\bf Proof } We have 
$\im(s) \subset \ker(d) = \im(\pi^*) \subset \ker(s)$
by the exactness of the Gysin sequence.
We have $w_{2i-1}=s(w_{2i})= \pi^*d(w_{2i}) \in 
\im(\pi^*)$ and as $d(w_{2i}^2)=0$, 
$w_{2i}^2\in \ker(d) = \im(\pi^*)$.
Hence we get the inclusion $A \subset \im(\pi^*)$,
where $A$ is the polynomial ring 
in variables $w_{2i-1}$ and $w_{2i}^2$. 
As already seen, $\im(s) \subset \im(\pi^*)$, so 
$\im(s) + A \subset \im(\pi^*)$. 
This completes the proof, as
$\ker(s) = \im(s) + A$ by Theorem \ref{y2k}. 
\hfill$\Box$

\medskip

\centerline{\bf Generators for the ring $H^*(BGO(2n))$}

\begin{lemma}\label{generators for H^*(BGO(2n))} 
The ring $H^*(BGO(2n))$ is generated 
by $\lambda,\, (a_{2i-1})_i,\, (b_{4i})_i,\,(d_T)_T
\in H^*(BGO(2n))$.
\end{lemma}

{\bf Proof } Let $S\subset H^*(BGO(2n))$ be the subring generated 
by the elements \\
$\lambda,\, (a_{2i-1})_i,\, (b_{4i})_i,\,(d_T)_T$.
By Theorem \ref{y2k}, the ring $B$
is generated by the elements $w_{2i-1} = \pi^*(a_{2i-1})$ and 
$w_{2i}^2 = \pi^*(b_{4i})$ where $1\le i\le n$, together
with the elements
$s(v_T)=\pi^*(d_T)$ where $T\subset \{ 1,\ldots,n\}$
with $|T| \ge 2$, which shows that 
$\pi^*(S) = B$.
As $1\in S$, we have $H^0(BGO(2n))\subset S$. 
Moreover, $H^1(BGO(2n))=\{ 0, a_1\} \subset S$.
We now proceed by induction. Suppose that
$H^j(BGO(2n))\subset S$ for all $j<i$. From $\pi^*(S) = B$,
and the fact that $\pi^*$ is a graded homomorphism, 
it follows that given $x\in H^i(BGO(2n))$ there exists
$x'\in S\cap H^i(BGO(2n))$ such that $\pi^*(x) = 
\pi^*(x')$. Hence $(x-x')\in \ker(\pi^*) = \im(\lambda)$,
so let $x = x' + \lambda y$ where $y\in H^{i-1}(BGO(2n))$.
By induction, $y\in S$, therefore $x\in S$. 
This proves the lemma. \hfill$\Box$

\begin{lemma}\label{algebraic independence} 
In the $k$-algebra $H^*(BGO(2n))$, the $n+1$ elements 
$\lambda$ and $b_{4i}$ (where $1\le i\le n$) are algebraically independent
over $k$. 
\end{lemma}

{\bf Proof } We recall the following commutative diagram
in which the top row is exact.
$$\begin{array}{ccccccc}
H^j(BGO(2n)) & \stackrel{\pi^*}{\to} & H^j(BO(2n)) 
& \stackrel{d}{\to} & H^{j-1}(BGO(2n))
& \stackrel{\lambda}{\to} & H^{j+1}(BGO(2n))\\
&&& \scriptstyle{s}\searrow ~~ &~~\dna\scriptstyle{\pi^*} & & \\
&&& & H^{j-1}(BO(2n)) & & 
\end{array}$$
Let $k[x_1,\ldots,x_n,y]$ be a polynomial ring in the
$n+1$ variables $x_1,\ldots,x_n,y$. Let
$f\in k[x_1,\ldots,x_n,y]$ be a non-constant polynomial of the
lowest possible total degree, with
$f(b_4,\ldots,b_{4n},\lambda) =0$.
Let $f = f_0 + yf_1$ where $f_0\in k[x_1,\ldots,x_n]$.
Then as $\pi^*(\lambda)=0$, we get
$$0= \pi^*f(b_4,\ldots,b_{4n},\lambda) = 
\pi^*f_0(b_4,\ldots,b_{4n}) = f_0(w_2^2,\ldots,w_{2n}^2)$$
As $w_{2i}^2$ are algebraically independent elements of
the $k$-algebra $H^*(BO(2n))$, the equality
$\pi^*(b_{4i})=w_{2i}^2$ implies that the 
$n$ elements $b_{4i}\in H^*(BGO(2n))$ 
are algebraically independent over $k$. 
It follows that $f_0 = 0 \in k[x_1,\ldots,x_n]$.
Hence we have $f= yf_1$, and so 
$\lambda f_1(b_4,\ldots,b_{4n},\lambda) = 
f(b_4,\ldots,b_{4n},\lambda) = 0$. Hence by exactness of the 
Gysin sequence, there exists $z \in H^*(BO(2n))$ with
$f_1(b_4,\ldots,b_{4n},\lambda) = dz$.
Applying $\pi^*$ to both sides, this gives
$\pi^*f_1(b_4,\ldots,b_{4n},\lambda) = s(z)$.
Now let $f_1 = f_2 + yf_3 \in 
k[x_1,\ldots,x_n,y]$, where $f_2 \in k[x_1,\ldots,x_n]$.
As $\pi^*f_1(b_4,\ldots,b_{4n},\lambda) = 
\pi^*f_2(b_4,\ldots,b_{4n})$, we get the equality
$\pi^*f_2(b_4,\ldots,b_{4n}) = s(z)$, that is,
$$f_2(w_2^2,\ldots ,w_{2n}^2) \in \im(s)$$
Now, from the formula $s = \sum_iu_i{\pa\over\pa v_i}$, 
it follows that $\im(s)$ is contained in the ideal 
generated by the $u_i = w_{2i-1}$. Hence the above
means that $f_2(w_2^2,\ldots ,w_{2n}^2) =0$.
Hence as before, $f_2 = 0 \in k[x_1,\ldots,x_n]$.
Hence we get $f = yf_1 = y^2f_3$. 

As $\lambda (\lambda f_3(b_4,\ldots,b_{4n},\lambda)) 
=  \lambda^2 f_3(b_4,\ldots,b_{4n},\lambda)=0$, 
by exactness of the Gysin sequence
there is some $q \in H^*(BO(2n))$ with
$d(q) = \lambda f_3(b_4,\ldots,b_{4n},\lambda)$.
Hence $s(q) = \pi^*d(q) = \pi^*(\lambda f_3(b_4,\ldots,b_{4n},\lambda))
=0$ as $\pi^*(\lambda)=0$. Hence $q \in \ker(s)$. 
The Lemma \ref{image of alpha} showed that $\ker(s) = \im(\pi^*)$,
hence $q = \pi^*(h)$ for some $h \in H^*(BGO(2n))$.
But then we have
$$\lambda f_3(b_4,\ldots,b_{4n},\lambda)  =
d(q) = d\pi^*(h) = 0$$
as $d\pi^*=0$ in the Gysin complex.
Hence the polynomial $g = yf_3 \in k[x_1,\ldots,x_n,y]$,
which is non-constant with degree less than that of $f = yg$,
has the property that $g(b_4,\ldots,b_{4n},\lambda) =0$.
This contradicts the choice of $f$, proving the lemma.
\hfill$\Box$

\rem\label{short.exact}
Let $k[b_4,\ldots,b_{4n},\lambda]$ denote the 
polynomial ring in the variables
$b_4,\ldots,b_{4n},\lambda$.
We have a short exact sequence of $k$-modules
$$0 \to k[b_4,\ldots,b_{4n},\lambda]
\stackrel{\lambda}{\to} H^*(BGO(2n))
\stackrel{\pi^*}{\to} B \to 0$$
where injectivity of $\lambda$ is by Lemma \ref{algebraic independence},
exactness in the middle is by exactness of the Gysin sequence, 
and surjectivity of $\pi^*$ is by Lemma \ref{image of alpha}.

\medskip

Now we state and prove the main result.

Consider the $2^n + n$
algebraically independent indeterminates
$\lambda$, $a_{2i-1}$ and $b_{4i}$ where $1\le i\le n$,
and $d_T$ where $T$ varies over subsets of $\{ 1,\ldots, n\}$ of
cardinality $|T|\ge 2$.
Let $k[\lambda, (a_{2i-1})_i, (b_{4i})_i, (d_T)_T ]$
be the polynomial rings in these variables. 

\begin{theorem}\label{cohomology} 
For any $n\ge 1$, the cohomology ring of $BGO(2n)$ with 
coefficients $k=\Z/(2)$ is isomorphic to the quotient 
$$H^*(BGO(2n))= 
{k[\lambda, (a_{2i-1})_i, (b_{4i})_i, (d_T)_T ]\over  
I}$$
where $I$ is the ideal generated by

{\bf (1) } the $n$ elements $\lambda a_{2i-1}$ for $1\le i\le n$, and

{\bf (2) } the $2^n -n -1$ elements $\lambda d_T$ 
where $T \subset \{ 1,\ldots,n\}$ is a subset of cardinality $\ge 2$, and

{\bf (3) } the $2^n -n(n-1)/2 -n-1$ elements of the form
$$\sum_{i\in T} a_{2i-1} \,d_{T-\{ i\}}$$
where $T \subset \{ 1,\ldots,n\}$ is a subset of cardinality $\ge 3$, and

{\bf (4) } the $n(n-1)/2$ elements of the form
$$(d_{\{ i,j \} })^2 + a_{2i-1}^2b_{4j} + a_{2j-1}^2b_{4i}$$
where $\{ i,j\}\subset \{ 1,\ldots,n\}$ is 
a subset of cardinality $2$, and finally

{\bf (5) } the $(2^n -n -1)^2 - n(n-1)/2$ elements (not necessarily distinct)
of the form
$$d_T\,d_U + \sum_{p\in T} \prod_{q \in T \cap U -\{ p\} }  
\, a_{2p-1} \, b_{4q} \, d_{(T - \{ p\}) \Delta U}$$
where $T\ne U$ in case both $T$ and $U$ have cardinality $2$,
and where $\Delta$ denotes the symmetric difference of sets 
$X \Delta Y = (X-Y)\cup (Y-X)$.
\end{theorem}

{\bf Proof } Let $R= k[\lambda, (a_{2i-1})_i, (b_{4i})_i, (d_T)_T ]$
denote the polynomial ring. 
Consider the homomorphism 
$R \to H^*(BGO(2n))$ which 
maps each of the variables in this polynomial ring
to the corresponding element of $H^*(BGO(2n))$. 
This map is surjective by Lemma \ref{generators for H^*(BGO(2n))}.

From  $\pi^*(\lambda)=0$ and from 
the description of $B$ in terms of
generators and relations given in Theorem \ref{y2k}, 
it follows that all the generators of the
ideal $I \subset R$,
which are listed above, map to $0$ under $\pi^* : H^*(BGO(2n)) \to B$.
Hence we get an induced surjective homomorphism
$$\varphi : R/I \to H^*(BGO(2n))$$

From its definition, we see that the ideal $I\subset R$ satisfies
$I \cap \lambda k[\lambda, b_4,\ldots,b_{4n}] =0$. Hence we get an inclusion
$\lambda k[\lambda, b_4,\ldots,b_{4n}] \hra R/I$.
Next, consider the map $R\to B$ which sends 
$\lambda \mapsto 0$, $a_{2i-1}\mapsto w_{2i-1}$,
$b_{4i}\mapsto w_{2i}^2$, and $d_T \mapsto c_T$. 
From  $\pi^*(\lambda)=0$ and from 
the description of $B$ in terms of
generators and relations given in Theorem \ref{y2k}, 
it again follows that this map is surjective, and 
all the generators of the ideal $I \subset R$
map to $0$, 
inducing a surjective homomorphism
$\psi : R/I \to B$. Hence we get a short exact sequence
$0 \to \lambda k[\lambda, b_4,\ldots,b_{4n}] \hra R/I
\stackrel{\psi}{\to} B \to 0$.

The above short exact sequence and the short exact sequence of
Remark \ref{short.exact} fit in the following commutative diagram.
$$\begin{array}{ccccccc}
0\to&\lambda k[\lambda,b_4,\ldots,b_{4n}]&\to&R/I&
\stackrel{\psi}{\to}&B&\to 0 \\
    & \|        &   &\scriptstyle{\varphi} \dna ~~ &   &\|    &      \\
0\to&\lambda k[\lambda,b_4,\ldots,b_{4n}]&\to&H^*(BGO(2n))&
\stackrel{\pi^*}{\to}&B&\to 0 
\end{array}$$
Hence the theorem follows by five lemma.
\hfill$\Box$

\example In particular, the small dimensional cohomology 
vector spaces of the
$BGO(2n)$ are as follows, in terms of linear bases over 
the coefficients $k=\Z/(2)$.
\begin{eqnarray*}
H^0(BGO(2n)) & = & <1> \\
H^1(BGO(2n)) & = & <a_1> \\
H^2(BGO(2n)) & = & <\lambda,\, a_1^2> \\
H^3(BGO(2n)) & = &  \left\{\begin{array}{ll}
                   <a_1^3> & \mbox{ for } n =1, \\
                   <a_1^3,\,a_3> & \mbox{ for } n\ge 2.
                   \end{array}\right.
                  \\
H^4(BGO(2n)) & = &  \left\{\begin{array}{ll}
                   <\lambda^2,\, a_1^4,\, b_4> 
                          & \mbox{ for } n =1, \\
                   <\lambda^2,\, a_1^4,\,a_1a_3,\, b_4>
                          & \mbox{ for } n\ge 2.
                    \end{array}\right.
                  \\
H^5(BGO(2n)) & = & \left\{\begin{array}{ll}
                   <a_1^5,\, a_1b_4>
                          & \mbox{ for } n =1, \\
                   <a_1^5,\, a_1^2a_3,\, a_1b_4,\, d_{\{1,2\}}>
                          & \mbox{ for } n =2, \\
                   <a_1^5,\, a_1^2a_3,\, a_1b_4,\, a_5, \, d_{\{1,2\}}>
                          & \mbox{ for } n\ge 3. 
                    \end{array}\right.
\end{eqnarray*}

\section{The map $H^*(BGL(2n))\to H^*(BGO(2n))$}

The cohomology ring $H^*(BGL(2n))$ is the polynomial ring
$\Z/(2)[\ov{c}_1, \ldots, \ov{c}_{2n}]$ where the $2n$ variables
$\ov{c_i}$ are the Chern classes mod $2$.
By definition of $b_{4r} \in H^*(BGO(2n))$ in terms of the  
universal triple $(E,L,b)$,
we have 
$\ov{c}_{2i}(E) = b_{4i}$ for $1\le i\le n$. 
Hence under the ring homomorphism
$H^*(BGL(2n))\to H^*(BGO(2n))$ induced by the inclusion 
$GO(2n)\hra GL(2n)$, we have 
$\ov{c}_{2i}\mapsto b_{4i}$. The following proposition gives
the images of the odd Chern classes $\ov{c}_{2i}$ in $H^*(BGO(2n))$,
which completes the determination of the ring homomorphism 
$H^*(BGL(2n))\to H^*(BGO(2n))$.

\begin{proposition}\label{odd Chern classes}
Consider the $n\times n$ matrix $A$ over $\Z/(2)[\lambda^2]$, 
with entries
$$A_{r,k} = {n-k \choose 2r - 2k} \lambda^{2r-2k}$$
This is a lower triangular matrix, with all diagonal entries equal to $1$,
so is invertible over the ring $\Z/(2)[\lambda^2]$.
Let $B = A^{-1}$ be its matrix inverse, which is again lower triangular,
with all diagonal entries $1$. For $1\le r\le n$, let the polynomial
$f_r(\lambda, b_4,\ldots, b_{4r-4})$ be defined by
$$f_r = {n \choose 2r-1} \lambda^{2r-2} 
+\sum_{1\le k\le r-1}{n-k \choose 2r-1-2k}\lambda^{2r-2-2k}
\left(\sum_{1\le j\le k}B_{k,j}(b_{4j} - {n \choose 2j} \lambda^{2j})\right)$$
By definition, $f_r(\lambda, b_4,\ldots, b_{4r-4})$ is
linear in  $b_4,\ldots,b_{4r-4}$, 
with the coefficient of $b_{4r-4}$ equal to the constant 
$n-r+1 \in \Z/(2)$. Then in the cohomology ring $H^*(BGO(2n))$, we have 
the identity
$$\ov{c}_{2r-1}(E) =  
a_{2r-1}^2  + \lambda \cdot f_r(\lambda, b_4,\ldots, b_{4r-4})$$
\end{proposition}

\proof We divide the proof into two steps, {\bf (a)} and  {\bf (b)}.

{\bf (a) } For each $1\le r\le n$, there exists a unique polynomial 
$g_r(\lambda, b_4,\ldots, b_{4r-4})$ such that
$\ov{c}_{2r-1}(E) =  
a_{2r-1}^2  + \lambda \cdot g_r(\lambda, b_4,\ldots, b_{4r-4})$.

{\bf (b) } The polynomials 
$g_r(\lambda, b_4,\ldots, b_{4r-4})$ are the polynomials 
$f_r(\lambda, b_4,\ldots, b_{4r-4})$ defined 
in the statement of the proposition.

{\bf Proof of (a) } The composite homomorphism
$O(2n)\hra GO(2n)\hra GL(2n)$
induces $H^*(BGL(2n))\to H^*(BGO(2n))\to H^*(BO(2n))$
under which $\ov{c}_i \mapsto w_i^2$.
As $a_{2i-1}^2\mapsto w_{2i-1}^2$ under $H^*(BGO(2n))\to H^*(BO(2n))$,
we must have 
$$\ov{c}_{2i-1}(E) = a_{2i-1}^2 + x_{4i-2}$$ 
where $x_{4i-2}$ lies in the kernel of 
$\pi^*: H^{4i-2}(BGO(2n))\to H^{4i-2}(BO(2n))$.
By exactness of Gysin, we have $x_{4i-2} = \lambda y_{4i-4}$
for some $y_{4i-4} \in H^{4i-4}(BGO(2n))$. 
Now, by the structure of the ring
$H^*(BGO(2n))$, we know that $\lambda$ annihilates 
the $a_{2j-1}$'s and the $d_T$'s. Hence we can replace 
$y_{4i-4}$ by a polynomial $g_i(\lambda, b_{4j})$ 
in $\lambda $ and the $b_{4j}$'s. As by Lemma \ref{algebraic independence}
the variables $\lambda$ and the $b_{4j}$ are algebraically independent, 
$g_i$ is unique. 
By degree considerations, the highest $b_{4j}$ that occurs 
in $g_i$ can be at most $b_{4i-4}$. This completes the proof of (a).

{\bf Proof of (b) } We first determine 
the polynomial $g_i(\lambda, b_4,\ldots, b_{4i-4})$ in the following
Example \ref{lambda and b} of a triple $\T = (E,L,b)$, which 
generalizes the Example \ref{eachL}. 

\example\label{lambda and b}
Let $\O(1)$ be the universal line bundle on $B(\C^*) = \PP^{\infty}_{\C}$.
Let $X = B(\C^*)\times \ldots \times B(\C^*)$ be the product of $n+1$ copies,
and let $p_i : X \to B(\C^*)$ be the projections for $0\le i \le n$.
Let $L = p_0^*(\O(1))$, and for $1\le i\le n$ let
$K_i = p_i^*(\O(1))$. On $X$, we get the nondegenerate triples
$\T_i = ((L\otimes K_i)\oplus K_i^{-1}, L, b_i)$ where $b_i$ is induced by
the canonical isomorphism 
$(L\otimes K_i) \otimes K_i^{-1} \stackrel{\sim}{\to} L$. 
Let $\T = (E,L,b)$ be the triple $\oplus \T_i$. 
We now write down the characteristic classes of $\T$. 
Let $\lambda = \ov{c}_1(L)$. Then by definition, 
$\lambda(\T)= \lambda$. 
As the odd cohomologies of $X$ are zero, the classes $a_{2i-1}(\T)$
and $d_T(\T)$ are zero, for all $1\le i \le n$ and for all
$T\subset \{ 1,\ldots, n\}$ with $|T| \ge 2$. 

Let $x_i = \ov{c}_1(K_i)$ for $1\le i\le n$, and let 
$s_k(\lambda x_i+ x_i^2)$ denote the $k$ th elementary symmetric
polynomial in the variables $\lambda x_1+ x_1^2,\ldots, \lambda x_n+ x_n^2$.
As $E$ is the direct sum of the $(L\otimes K_i) \oplus K_i^{-1}$,
its mod $2$ Chern classes $\ov{c}_i$ are given by 
\begin{eqnarray*}
\ov{c}_{2r}(E) & = & {n \choose 2r}\lambda^{2r} +
                 \sum_{1\le k\le r}{n-k \choose 2r-2k} \lambda^{2r-2k}
                 s_k(\lambda x_i+ x_i^2) \\
\ov{c}_{2r-1}(E) & = & {n \choose 2r-1} \lambda^{2r-1} +
\sum_{1\le k\le r-1}{n-k \choose 2r-1-2k} \lambda^{2r-1-2k}
                 s_k(\lambda x_i+ x_i^2) 
\end{eqnarray*}
As $b_{4r}(\T) = \ov{c}_{2r}(E)$, we have the following
equations for $1\le r\le n$.
$$b_{4r} -{n \choose 2r}\lambda^{2r}  
= \sum_{1\le k\le r}{n-k \choose 2r-2k} \lambda^{2r-2k} 
s_k(\lambda x_i + x_i^2)$$
This is a system of linear equations with coefficients in $\Z/(2)[\lambda^2]$, 
for the $b_{4r}-{n \choose 2r} \lambda^{2r}$ in terms of the 
$s_k(\lambda x_i + x_i^2)$. 
It is given by the $n\times n$ matrix $A$ over $\Z/(2)[\lambda^2]$, 
with entries
$$A_{r,k} = {n-k \choose 2r - 2k} \lambda^{2r-2k}$$
This is a lower triangular matrix, with all diagonal entries equal to $1$,
so is invertible over the ring $\Z/(2)[\lambda^2]$.
Let $B = A^{-1}$ be its matrix inverse, which is again lower triangular,
with diagonal entries $1$. Hence we get
$$s_r(\lambda x_i + x_i^2)  = \sum_{1\le k\le r}B_{r,k}
(b_{4k} - {n \choose 2k} \lambda^{2k})$$
Substituting this in the equation for 
$\ov{c}_{2r-1}(E)$, we get equations 
$$ \ov{c}_{2r-1}(E) = \lambda \cdot f_r(\lambda, b_4,\ldots, b_{4r-4})$$
where 
$$f_r = {n \choose 2r-1}\lambda^{2r-2} +
\sum_{1\le k\le r-1}{n-k \choose 2r-1-2k}
             \lambda^{2r-2-2k}\left(
\sum_{1\le j\le k}B_{k,j}
(b_{4j} - {n \choose 2j} \lambda^{2j})\right)$$

\medskip

{\bf Proof of (b) continued } In the Example \ref{lambda and b},
we have the desired equality $g_i = f_i$.
Note that the cohomology ring $H^*(X)$ is the polynomial ring
$\Z/(2)[\lambda, x_1,\ldots,x_n]$, in which 
the $n+1$ elements $\lambda, b_4,\ldots, b_{4n}$ 
are algebraically independent, where 
$b_{4r} = {n \choose 2r}\lambda^{2r} +
\sum_{1\le k\le r}{n-k \choose 2r-2k} \lambda^{2r-2k} 
s_k(\lambda x_i+ x_i^2)$. Hence as
$g_i = f_i$ in this example, we get $g_i = f_i$ universally.

This completes the proof of Proposition \ref{odd Chern classes}.
\hfill$\Box$

\medskip

\example The above proposition in particular gives in $H^*(BGO(2n))$
the identities
\begin{eqnarray*}
\ov{c}_1 & = & a_1^2 + n \lambda ~~\mbox{for all } n \ge 1, ~\mbox{and}\\
\ov{c}_3 & = & a_3^2 + {n(n-1)(2n-1) \over 6} \lambda^3 + (n-1)\lambda b_4 
                    ~~\mbox{for all } n \ge 2.
\end{eqnarray*}

\medskip

{\bf Acknowledgement } Y. I. Holla thanks ICTP Trieste, and 
N. Nitsure thanks the University of Essen, for support while
this work was being completed.

\section*{References} \addcontentsline{toc}{section}{References}

[H-N] Holla, Y. I. and Nitsure, N. : Characteristic Classes for
$BGO(2n)$. (Earlier version of this paper, without the last section).
To appear in Asian J. Mathematics.

[M] Matsumura, H. :  {\sl Commutative ring theory}, Cambridge 
Studies in Advanced Mathematics {\bf 8}, Cambridge Univ. Press, 1986.

\bigskip

\bigskip

\bigskip

{\it
Address : School of Mathematics, Tata Institute of Fundamental
Research, Homi Bhabha Road, Mumbai 400 005, India.

e-mails: yogi@math.tifr.res.in, nitsure@math.tifr.res.in 
}

\bigskip

\centerline{22-III-2000, revised on 22-V-2000}

\end{document}